\documentclass[twocolumn]{autart} 
\usepackage[utf8]{inputenc}
\usepackage{amsfonts, amsmath, amssymb}
\usepackage[table]{xcolor} 
\usepackage[shortlabels]{enumitem}  
\usepackage[font=footnotesize]{caption, subcaption}  
\usepackage[T1]{fontenc} 
\usepackage{graphicx} 
\usepackage{tabularx}  
\usepackage{cite} 

\usepackage{algorithm}
\usepackage{algpseudocode}
\usepackage{algorithmicx}  

\algnewcommand{\AlgAnd}{\textbf{ and }}

\setlist{noitemsep} 
\setlist[enumerate,1]{label = (\roman*)} 

\newcommand\restr[2]{\ensuremath{\left.#1\right|_{#2}}}
 
\newcommand{\edit}[1]{{\color{black}#1}} 
\DeclareMathOperator*{\argmin}{arg\,min}

\begin{document}
\begin{frontmatter}

\title{Auction algorithm sensitivity for multi-robot task allocation}
\thanks{This project was partially supported by TAS-DCRC.}

\author[UNSW]{Katie Clinch}\ead{k.clinch@unsw.edu.au},
\author[EPFL]{Tony A.\ Wood}\ead{tony.wood.@epfl.ch}, 
\author[Melb]{Chris Manzie}\ead{manziec@unimelb.edu.au}
\address[UNSW]{School of Computer Science and Engineering,  University of New South Wales, Sydney, Australia.}
\address[EPFL]{SYCAMORE Lab, \'Ecole Polytechnique F\`ed\`erale de Lausanne (EPFL), Lausanne, Switzerland.}
\address[Melb]{Department of Electrical and Electronic Engineering, University of Melbourne, Parkville, Victoria, 3010, Australia.}

\begin{keyword}
 Sensitivity analysis, robustness, auction algorithm, minimum spanning tree, traveling salesman problem, task assignment, approximation algorithm.
\end{keyword}

\begin{abstract}
We consider the problem of finding a low-cost allocation and ordering of tasks between a team of robots in a $d$-dimensional, uncertain, landscape, and the sensitivity of this solution to changes in the cost function.

Various algorithms have been shown to give a 2-approximation to the \textsc{MinSum} allocation problem \cite{L+_AuctionBasedMRR}. By analysing such an auction algorithm, we obtain  intervals on each cost, such that any fluctuation of the costs within these intervals will result in the auction algorithm outputting the same solution.
\end{abstract}
\end{frontmatter}

\section{Introduction} \label{sec:Intro}

The question of how best to allocate tasks at specified locations between a team of mobile robots is a key problem in the study of \emph{multi-robot systems}. 
Questions of this type are commonly referred to as \emph{multi-robot task allocation} (MRTA) problems \cite{GM_taxonomy, KHE_MRTAOverview}
and are related to the \emph{multiple travelling salesman problem} \cite{B_Survey_mTSP}, 
\emph{vehicle routing problem} \cite{TV_survey_VRP} 
and \emph{generalised assignment problem} \cite{Pentico_survey_AP}. 
\edit{
This NP-hard question has been of interest for decades, and has inspired a wealth of literature 
for finding good, sub-optimal, solutions 
\cite{AB_MILP, 
KNGG_MonteCarlo, 
B_auctionAlgorithm,Z+_MRTAMarketEconomy, 
GCY_genetic, 
WGL_AntColony, 
MM_SimulatedAnnealing}. 

A second question is: given an MRTA problem and an algorithm which outputs a good, sub-optimal, solution to it, how much can the inputs change without effecting the output of the algorithm? 
Unlike the first question, this second question (which we call the \emph{sensitivity problem}) has received little attention. 

Knowledge of the output's sensitivity is important when the online conditions may vary from the input data used in evaluating the offline solution to the MRTA problem. This can occur if the task locations are estimated; the terrain data is uncertain; or there are moving obstacles that necessitate local collision avoidance and subsequently longer path lengths. Knowing how much the incurred costs can vary before the MRTA solution changes can avoid unnecessary online replanning. Additionally, understanding which parts of the algorithm's solution are most  sensitive to inaccuracies can facilitate targeted changes in resourcing to provide more robust solutions. For example, selected robots may be supplied with extra fuel or enhanced communication abilities for specific segments. 

}

\edit{In this paper, we build on recent advances in sensitivity analysis to construct a theory to address the sensitivity problem quickly for a simple MRTA algorithm. 
We first show that, despite its simplicity, our proposed algorithm  \textsc{Assign}} obtains the same guarantee on sub-optimality as the best MRTA algorithms considered in \cite{K+_PowerOfSSIAuctions,L+_AuctionBasedMRR}. 

\edit{We then address the sensitivity problem, by considering the effect of input errors on the cost function for \textsc{Assign}. Namely, we show how to construct an interval for each cost, so that any fluctuation of costs within their corresponding intervals will not change the output of $\textsc{Assign}$.

There are in fact many families of intervals which solve the sensitivity question for \textsc{Assign}. Based on the idea that the best choice of interval family will first maximise the smallest margin of error, then the second smallest etc., we finally prove our main result: how to obtain the \emph{lexicographic best} family of intervals for the sensitivity problem.} 


\subsection{Paper outline}

We collate key mathematical ideas in Section \ref{sec:Prelims}, before formally stating the \textsc{MinSum} problem in Section \ref{sec:Problem}. In Section \ref{sec:Auction}, we describe the algorithm \textsc{Assign} and show that it outputs a solution whose cost is at most twice the optimal. Section \ref{sec:Sensitivity} contains our main contributions. We analyse the sensitivity of \textsc{Assign} to both a single cost changing (Subsection \ref{subsec:OneEdge}), and all costs changing simultaneously (Subsection \ref{subsec:AllEdges}), and 
finally show how to obtain the lexicographic best family of intervals (Subsection \ref{subsec:AllEdges_Best}). 
Section \ref{sec:Simulations} puts this theory into practice, by providing results from computer simulations. Finally, in Section \ref{sec:FurtherWork}, we suggest future research directions.

\subsection{Related work} 

MRTA problems fit under the broader category of \emph{assignment problems}, which have been studied since Kuhn's seminal paper \cite{K_hungarian} introduced the Hungarian algorithm for the classic assignment problem. \emph{Sensitivity analysis}, which determines how much inputs can fluctuate before an algorithm outputs a different solution, has also been studied since the 1950's \cite{WW_SensitivitySurvey}. 

Perhaps surprisingly, it is only relatively recently that these two theories have been combined: Mills-Tettey, Stentz and Dias consider the sensitivity of the Hungarian algorithm to changes in the cost function \cite{MT+_SensitivityAP}.
Similar techniques have been used to analyse the cost-sensitivity of algorithms for the 
bottleneck assignment problem \cite{M+_SensitivityBottleneckAP}, 
linear assignment problem \cite{LW_SensitivityLAP, M+_SensitivityLAP} 
and the MRTA problem where each robot is assigned at most one task \cite{LS_MRTASensitivity_Hungarian}.
To the best of our knowledge, this is the first paper to apply sensitivity analysis to an MRTA algorithm which assigns multiple tasks to each robot. 

\section{Preliminaries}\label{sec:Prelims}

\subsection{Notation} \label{subsec:Notation}

We use $\{\cdot\}$ to denote \emph{sets}, and extend vector notation by allowing $(\cdot)$ to denote ordered \emph{multisets}, also known as \emph{lists}.
We set $\mathbb{R}_+ = \{x\in \mathbb{R}: x\geq 0 \}\cup\{\infty\}$.

Our problem setting consists of a set of \emph{robots} $R=\{r_i\}_{i=1}^{m}$ and \emph{tasks} $T=\{t_j\}_{j=1}^n$ in a 
\emph{landscape} 
$L\subseteq \mathbb{R}^d$ for some positive integer \emph{dimension} $d$. Each robot $r_i$ and task $t_j$ has a corresponding \emph{initial position} $p(r_i)\in L$ and \emph{position} $p(t_j)\in L$ respectively. 
%
For any distinct pair $p_1,p_2\in L$ there is a non-negative \emph{traversal cost} $c(p_1,p_2)$ of traversing the landscape from $p_1$ to $p_2$. We assume this cost is identical for all robots. 
To simplify notation, we set $c(s_i, s_j) = c(p(s_i), p(s_j))$ for all distinct $s_i,s_j\in R\sqcup T$. 

Depending on the setting, there may also be a non-negative \emph{execution cost} $c(s,s)$ corresponding to  executing a task $s\in T$ or booting-up a robot $s\in R$. 

The mapping $c:(R\sqcup T)\times (R\sqcup T) \to \mathbb{R}_+$ that results from the traversal and execution costs is called the \emph{cost function}. 
An \emph{instance} of an MRTA problem with robots $R$, tasks $T$, landscape $L$, positions $p$,  and cost function $c$ is denoted  $\mathcal{M}=(R,T,L,p,c)$, and has $|R\sqcup T| \geq 3$. 

\begin{rem}
\edit{In Section \ref{subsec:MetricCost}, we restrict to problem settings with no execution costs (i.e. $c(s,s) = 0$ for all $s\in R\sqcup T$) for the remainder of the paper. 
Additionally, the algorithms we later describe do not use inter-robot costs ($c(r_1,r_2)$ where $r_1,r_2\in R$ are distinct).
So 
to implement the results in this paper, it suffices to know 
$c(s_1,s_2)$ for all $s_1\in R\cup T$ and $s_2\in T\setminus\{s_2\}$.}
\end{rem}

\subsection{Metric instances} \label{subsec:MetricTSP}

The cost function $c$ of an instance $\mathcal{M} = (R,T,L,p,c)$ is \emph{metric} if for all $s_1,s_2,s_3\in R\sqcup T$,
\begin{enumerate}[label = (M\arabic*)]
    \item $c(s_1,s_2)\geq 0$ ($c$ is non-negative), \label{defn:metric_positive}
    \item $c(s_1,s_2) = 0$ if and only if $s_1=s_2$, \label{defn:metric_zero}
    \item $c(s_1,s_2) = c(s_2,s_1)$ ($c$ is \emph{symmetric}), and \label{defn:metric_symm}
    \item $c(s_1,s_3) \leq c(s_1,s_2) + c(s_2,s_3)$ ($c$ satisfies the \emph{triangle inequality}). \label{defn:metric_triangle}
\end{enumerate}
If \ref{defn:metric_zero} is weakened to the following, then $c$ is \emph{pseudometric}.
\begin{enumerate}
    \item[(M2')] $c(s,s) = 0$ for all $s\in R\sqcup T$. \label{defn:pseudo_zero}
\end{enumerate}
%
We say an instance $\mathcal{M} = (R,T,L,p,c)$ is \emph{(pseudo)metric} when $c$ is (pseudo)metric.

\subsection{Optimization and approximation} \label{subsec:Optimization} 


\edit{

An \emph{MRTA problem} $\textsc{Prob}$ on an instance $\mathcal{M} = (R,T,L,p,c)$ minimises a non-negative function $C$ (derived from $c$) on $\mathcal{M}$ subject to certain constraints. An \emph{optimal solution} to $\textsc{Prob}$ on $\mathcal{M}$ is denoted $\textsc{Prob}(\mathcal{M})$ and has \emph{cost} $C(\textsc{Prob}(\mathcal{M}))$.

A polynomial-time algorithm $\textsc{Alg}$ is said to be a \emph{$k$-approximation} for $\textsc{Prob}$ (or equivalently have \emph{performance ratio} $k$) if for all instances $\mathcal{M}$, the \emph{solution} $\textsc{Alg}(\mathcal{M})$ satisfies the constraints for $\textsc{Prob}$ and has $C(\textsc{Alg}(\mathcal{M})) \leq k\cdot C(\textsc{Prob}(\mathcal{M}))$.

Several optimization problems which are NP-hard to solve in general, are comparatively easy to approximate when the instance is (pseudo)metric (see Section \ref{subsec:MetricCost}).  
We shall see that this is the case for the metric \textsc{MinSum} problem (defined in Section \ref{sec:Problem}), which has 2-approximation algorithm \textsc{Assign} (described in Section \ref{sec:Auction}).
}




\subsection{Graphs}\label{subsec:Graphs}

A \emph{graph} $G=(V;E)$ consists of a finite set of \emph{vertices} $V(G)=V$, and \emph{edges} $E(G)=E\subseteq\{\{u,v\}: u,v\in V\text{ are distinct}\}$. 
Given a vertex $v\in V$, the \emph{neighbourhood of $v$}, $N_G(v)$, is the set of all vertices $u\in V$ such that $\{u,v\}\in E$. 
A mapping $l: X \to \mathbb{R}_+$ is called respectively a \emph{vertex-labelling}, \emph{pair-labelling}, or \emph{edge-labelling} of $G$ when $X=V$, $X= V\times V$, or $X=E$ respectively. 

A list $P=(v_0,v_1,\ldots,v_n)$ of vertices $v_i$ of a graph $G$ is a \emph{route} (more commonly called a \emph{walk}) if $\{v_{i},v_{i+1}\}\in E(G)$ for all $i\in\{0, \dots, n-1\}$; in which case $P$ has edge set $E(P) = \{\{v_{i},v_{i+1}\}: 0\leq i\leq n-1\}$. A route with no repeated edges is a \emph{trail}, and a trail with no repeated vertices is a \emph{path}.
A trail $(v_0,\dots, v_n)$ whose only repeated vertex is $v_0 = v_n$ is called a \emph{cycle}. 
A graph containing no cycles is a \emph{forest}. If $G$ is a forest and for every $u,v\in V(G)$ there is a path in $G$ containing both $u$ and $v$, then $G$ is a \emph{tree}.

We model an MRTA instance $\mathcal{M} = (R,T,L,p,c)$ with a \emph{robot-task graph} $G(\mathcal{M})=(R,T; E;c)$, where $V(G(\mathcal{M})) = R\sqcup T$, $E(G(\mathcal{M})) = E = \{\{s_1,s_2\} : s_1,s_2\in R\sqcup T \text{ are distinct}\}$, and $c$ is the pair-labelling given by the cost function. 
Observe that when $c$ is (pseudo)metric we can simplify this model to the edge-labelled graph $G(\mathcal{M})=(R,T;E;\restr{c}{E})$.
\section{The \textsc{MinSum} problem} \label{sec:Problem}

Let $\mathcal{M} = (R,T,L,p,c)$ be an MRTA instance, and $P=(v_i)_{i=0}^n$ be a route in the corresponding robot-task graph $G$.  We define the \emph{cost of $P$} as
\begin{equation}\label{eqn:PathCostDefn}
\edit{C}(P) = \sum_{i=0}^{n} c(v_i,v_i) + \sum_{i=0}^{n-1} c(v_i, v_{i+1})
\end{equation}
when $n\geq 1$, and $\edit{C}(P)=0$ otherwise. 
Observe that when $\mathcal{M}$ is (pseudo)metric, this simplifies to $\edit{C}(P) = \sum_{i=0}^{n-1} c(v_i, v_{i+1})$. The route $P=(v_i)_{i=0}^n$ is a \emph{robot-route} if $v_0\in R$ and $v_j\in T$ for all $j\geq 1$.
%
A \emph{solution} to $\mathcal{M}$, also called a \emph{plan}, is a family $\mathcal{P} = \{P_1, \ldots, P_{|R|}\}$ of \emph{robot-routes} in $G$ such that $\mathcal{P}$ partitions $V(G)$. 
For the $\textsc{MinSum}$ problem, 
the \emph{cost of a plan} $\mathcal{P} = \{P_i\}_{i=1}^{|R|}$ 
is defined as
$
\edit{C}(\mathcal{P}) = \sum_{P_i\in \mathcal{P}} \edit{C}(P_i);
$
and $\textsc{MinSum}(\mathcal{M}) = \mathcal{P}^*$ for some plan 
$\mathcal{P}^* = \arg \min_\mathcal{P} \edit{C}(\mathcal{P})$.

In Section \ref{sec:Auction}, we  describe an approximation algorithm $\textsc{Assign}$ for $\textsc{MinSum}$ on metric instances $\mathcal{M}$, and then analyse the sensitivity of $\textsc{Assign}(\mathcal{M})$ in Section \ref{sec:Sensitivity}. Before this, we \edit{use the remainder of this section to justify} 
our choice to restrict to metric cost functions. This is a realistic assumption for many MRTA problems, but as we shall see, it is also justified by the current state of the theory.

Note that we can define other optimization problems on $\mathcal{M}$ such as $\textsc{MinMax}$ which minimises the cost of the most expensive route; or $\textsc{MinAvg}$ which minimises the average route cost. These formulations are also NP-hard. 

\subsection{Connection to the traveling salesman problem}\label{subsec:MetricCost}

If $\mathcal{M}$ is an instance with exactly one robot, then the problems $\textsc{MinSum}$, $\textsc{MinMax}$ and $\textsc{MinAvg}$ on $\mathcal{M}$ are equivalent to the infamous travelling salesman problem (TSP). 
Thus for instances $\mathcal{M}$ with multiple robots, these MRTA problems are extensions of TSP. As such, any properties which make TSP challenging are likely to make these MRTA problems more challenging still.

Christofides's \cite{Christofides} long-standing \cite{KKGS_2021_ImprovedMetricTSP} 
result provides a $\frac{3}{2}$-approximation algorithm for TSP when the cost function $c$ is (pseudo)metric. In comparison, if $c$ is allowed to fail \ref{defn:metric_triangle} then it is NP-complete to find a $k$-approximation algorithm for TSP, where $k\geq 1$ is any constant \cite{SG_PCompleteProblems}. The setting where condition \ref{defn:metric_symm} is relaxed is called asymmetric TSP (ATSP), and was only recently shown to have a constant bound approximation algorithm \cite{STV_AsymmetricTSP}. The current best known performance ratio for ATSP is 22 (see \cite{TV_ImprovedATSP}), which is too large to be of practical interest. 

Cost functions which arise from real-world data typically satisfy \ref{defn:metric_positive}. 
Together with the above results, this means it is sensible to restrict our analysis of \textsc{MinSum} to instances $\mathcal{M}=(R,T,L,p,c)$ where $c$ satisfies \ref{defn:metric_positive}, \ref{defn:metric_symm} and \ref{defn:metric_triangle}. 
Lemma \ref{lem:CanAvoidVertexLabelling}, below, 
shows that such an instance $\mathcal{M}$ whose cost function fails \ref{defn:metric_zero} 
can be restated in terms of either a metric cost function (in which case we can apply the methods developed in this paper) or an asymmetric cost function (in which case no good approximation algorithm is known). Hence justifying our choice to restrict to metric instances for the remainder of the paper.

\begin{lem}\label{lem:CanAvoidVertexLabelling}
Let $\mathcal{M} = (R,T,L,p,c)$ be an MRTA instance where $c$ satisfies \ref{defn:metric_positive}, \ref{defn:metric_symm} and \ref{defn:metric_triangle}, and let $G$ 
be the robot-task graph for $\mathcal{M}$. 
Then there exists a graph $G'$ with pair-labelling $c'$ such that
\begin{enumerate}
    \item $V(G')$ is a partition of $V(G)$, \label{part:metric_partition}
    \item $c'$ satisfies  \ref{defn:metric_positive} and \ref{defn:metric_zero}, and \label{part:metric_M12}
    \item  $\edit{C}'(P') = \edit{C}(P)$ for all robot-routes $P$ of $G$ (\edit{where $C$ and $C'$ are derived from $c$ and $c'$ respectively by the definition given in equation \eqref{eqn:PathCostDefn}, and } $P'$ is obtained from $P$ by replacing each $v\in P$ with the corresponding set $V_i\in V(G')$ for which $v\in V_i$).\label{part:metric_routes}
\end{enumerate} 
%
Further, $c'$ is metric if and only if there exists some $\lambda\in \mathbb{R}_+$ such that $c(t,t)=\lambda$ for all $t\in T$.
\end{lem}

\begin{pf}
\edit{If $c$ satisfies \ref{defn:metric_zero}, then set $G'=G$ and $c'=c$ and we are done. So suppose not. Since $c$ satisfies \ref{defn:metric_positive}, there are exactly two ways it can fail \ref{defn:metric_zero}: either there exist distinct $x,y\in V(G)$ such that $c(x,y) = 0$, or there is some $x\in V(G)$ for which $c(x,x)>0$. To ensure $c'$ satisfies \ref{defn:metric_zero}, we construct $c'$ from $c$ in two steps, each of which corrects for one of these two respective issues.}
%
First we form the partition $\mathcal{V} = V(G')$ of $V(G)$ and show that $c$ induces a function $c_\mathcal{V}$ on $\mathcal{V}\times \mathcal{V}$. Then we define $c'$ in terms of $c_\mathcal{V}$.

Observe that if $x,y\in V(G)$ are distinct with $c(x,y) = 0$, then \ref{defn:metric_positive}, \ref{defn:metric_symm} and \ref{defn:metric_triangle} imply 
$c(x,x) = c(y,y) = 0$. Similarly if $x,y,z\in V(G)$ are distinct with $c(x,y) = c(y,z) = 0$ then \ref{defn:metric_positive} and \ref{defn:metric_triangle} give $c(x,z) = 0$. Thus we take  $\mathcal{V}$ to be the partition of $V(G)$ into maximal sets $V_i$ such that either $|V_i| =1$, or $c(x,y)=0$ for all 
$x,y\in V_i$.
Property \ref{defn:metric_triangle} implies that for any distinct $V_1,V_2\in\mathcal{V}$ and vertices $v_1,v_1'\in V_1$ and $v_2\in V_2$, we have $c(v_1,v_2) = c(v_1',v_2)$. 
And so we can define a function $c_\mathcal{V} : \mathcal{V} \times \mathcal{V}\to \mathbb{R}_+$ by $c_\mathcal{V}(V_1,V_2) = c(v_1,v_2)$ for all $V_1,V_2\in\mathcal{V}$, $v_1\in V_1$ and $v_2\in V_2$. 

Let $G'$ be the complete graph with $V(G') = \mathcal{V}$. 
We define the pair-labelling $c'$ on $G'$ by setting
$
c'(V_1,V_1) = 0
$
for all $V_1 \in V(G')$, and
\begin{multline}\label{eqn:c'Definition}
c'(V_1,V_2) =  \alpha(V_1,V_2) \cdot c_\mathcal{V}(V_1,V_1) + c_\mathcal{V}(V_1,V_2) \\
+ c_\mathcal{V}(V_2,V_2) 
\end{multline}
for all distinct $V_1,V_2\in V(G)$; where $\alpha(V_1,V_2) = 0$ if $V_i\cap T\neq\emptyset$ for both $i\in\{1,2\}$, and $\alpha(V_1,V_2) = 1$ otherwise.
By construction, $c'$ clearly satisfies \ref{defn:metric_positive} and \ref{defn:metric_zero}. Thus $G'$ and $c'$ satisfy properties \ref{part:metric_partition} and \ref{part:metric_M12} respectively. We next consider property \ref{part:metric_routes}.

Let $P=(v_0,\dots, v_n)$ be a robot-route in $G$. For each $v_i$ in $P$, let $V_i$ denote the set in $\mathcal{V}$ containing $v_i$. Then $P'=(V_0,\dots, V_n)$ and for all $i\geq 1$ we have $v_i\in V_i\cap T$. Thus, by the definitions of $c'$, $c_\mathcal{V}$, and the cost of a route given in equation \eqref{eqn:PathCostDefn}, we have
\begin{align*}
\edit{C}'(P') &= \sum_{i=0}^n c'(V_i, V_i) + \sum_{i=0}^{n-1} c'(V_i, V_{i+1}) \\
  &= c(v_0,v_0) + \sum_{i=0}^{n-1} c(v_i,v_{i+1}) + c(v_{i+1},v_{i+1}) = \edit{C}(P).
\end{align*}
%
To prove the final equivalence statement in the lemma, first suppose $c'$ is metric.
Then 
equation \eqref{eqn:c'Definition}, the definition of $c_\mathcal{V}$, and the fact $c$ satisfies \ref{defn:metric_symm} 
imply there is some $\lambda\in\mathbb{R}_+$ such that $c(t,t) = \lambda$ for all $t\in T$. 

Conversely, suppose $c(t,t) = \lambda$ for all $t\in T$, and let $V_1,V_2,V_3\in\mathcal{V}$. If $V_1=V_2$ or $\alpha(V_1,V_2) = 1$, then trivially $c'(V_1,V_2) = c'(V_2,V_1)$. 
Otherwise there exist tasks $t_1\in V_1$ and $t_2\in V_2$, giving $c'(V_1,V_2) = c(t_1,t_2) + \lambda$; and so $c'$ inherits property \ref{defn:metric_symm} from $c$. 
It only remains to show that $c'$ satisfies \ref{defn:metric_triangle}. Let $v_i\in V_i$ for all $i\in\{1,2,3\}$. Since $c$ satisfies \ref{defn:metric_triangle}, we have
\begin{multline} \label{eqn:c'M4}
c'(V_1,V_2) + c'(V_2,V_3) \geq [\alpha(V_1,V_2)\cdot c(v_1,v_1) + c(v_2,v_2)] \\
+ c(v_1,v_3) + c(v_3,v_3).    
\end{multline}
If $\alpha(V_1,V_2)\geq \alpha(V_1,V_3)$, then the right hand side of \eqref{eqn:c'M4} is at least $c'(V_1,V_3)$ and we are done. Otherwise, $V_1$ and $V_2$ both contain tasks, and so $c(v_2,v_2) = \lambda = c(v_1,v_1)$, and we are again done.\hfill$\square$
\end{pf}


\begin{rem}
Lemma \ref{lem:CanAvoidVertexLabelling} implies two things. Firstly, that our subsequent results for metric MRTA instances also apply to instances where tasks have non-zero (but identical) execution costs, and robots have non-zero (and possibly distinct) booting-up costs.

And secondly, that if an MRTA instance $\mathcal{M}$ contains tasks with different non-zero execution costs, then $\textsc{MinSum}(\mathcal{M})$ is an extension of ATSP. So our restriction to metric cost functions covers exactly those MRTA instances for which it is currently practical to implement approximation algorithms.
\end{rem}

The following well-known fact shows that our results can also be applied to MRTA instances which balance minimising various objectives (e.g.\ time, battery depletion, accumulated damage), or to explore the Pareto optimal of multiple objectives, so long as each individual objective defines a metric.

\begin{fact}
Let $\{c_i\}_{i=1}^n$ be a finite set of (pseudo)metrics defined on a set $X$, and let $(\alpha_i)_{i=1}^{n}\in\mathbb{R}_+^n$. 
Suppose $c(x,y) = \sum_{i=1}^n \alpha_i c_i(x,y)$ for all $x,y\in X$.
Then $c$ is a (pseudo)metric.
\end{fact}
\section{The task allocation algorithm} \label{sec:Auction}


In this section, we describe the algorithm \textsc{Assign} (Algorithm \ref{alg:Planner}) 
which outputs a plan $\mathcal{P}$ in two steps. First a sequential single-item auction greedily allocates tasks to robots (\textsc{Auction}), and second a depth-first search converts each robot's allocation to a robot-route ($\textsc{DFShortcut}$). These steps will be described in detail in Subsections \ref{subsec:Auction} and \ref{subsec:DoubleTreeAlgo} respectively.
We then show that \textsc{Assign} provides a 2-approximation for metric \textsc{MinSum} (see Subsection \ref{subsec:2Approx}). 
\textsc{Assign} can be thought of as a multi-robot extension of the classic double-tree method for metric TSP \cite{RSL_SeveralHeuristicsForTSP}.

\begin{algorithm}[t]
    \caption{\textsc{Assign}}
    \textbf{input:} 
    $\mathcal{M} = (R,T,L,p,c)$\\
    \textbf{output:} plan  $\mathcal{P} = \{P_i: P_i \text{ is a robot-route for } r_i\in R \text{ and } \bigsqcup_{i} V(P_i) = R\sqcup T\}$
\begin{algorithmic}[1]
    \State $(W,a) \gets \textsc{Auction}(\mathcal{M})$ 
    \State $\mathcal{P}\gets\textsc{DFShortcut}(R,T,W, a)$ 
    \State \Return $\mathcal{P}$ 
\end{algorithmic}\label{alg:Planner}
\end{algorithm}

In Section \ref{sec:Sensitivity}, we show that \textsc{AuctionSensitivity} (Algorithm \ref{alg:Sensitivity}) can be used calculate the sensitivity of $\textsc{Auction}(\mathcal{M})$ to changes in the cost function. This in turn will give sufficient conditions for the output of $\textsc{Assign}$ to remain unchanged, despite changes in the cost function. We formalise these ideas in Section \ref{sec:Sensitivity}. The subroutines $\textsc{Initialiser}$ and $\textsc{ErrorIntervals}$ of Algorithm \ref{alg:Sensitivity} are described in Subsections \ref{subsec:AllEdges_Best} and \ref{subsec:AllEdges} respectively.

\begin{algorithm}[t]
    \caption{\textsc{AuctionSensitivity}}
    \textbf{input:} 
    $\mathcal{M} = (R,T,L,p,c)$, $\textsc{Auction}(\mathcal{M})=(W,a)$\\
    \textbf{output:} function $I$ 
    with $I(e) = (\underline{\Delta}(e),\overline{\Delta}(e)]$ for all edges $e$ of the robot-task graph
\begin{algorithmic}[1]
    \State $E \gets \{\{s_1,s_2\}: s_1,s_2\in R\sqcup T \text{ and } s_1\neq s_2 \}$ 
    \State $G \gets (R,T;E;c)$ 
    \State $I_0 \gets \textsc{Initialiser}(G,W,a)$ 
    \State $I\gets \textsc{ErrorIntervals}(G,W,a,I_0)$ 
    \State \Return $I$
\end{algorithmic}\label{alg:Sensitivity}
\end{algorithm}


\subsection{Step 1: \textsc{Auction}}\label{subsec:Auction}

\begin{algorithm}[b]
    \caption{\textsc{Auction}}
    
    \textbf{input:} 
    $\mathcal{M} = (R,T,L,p,c)$ 
    \\
    \textbf{output:} 
    list of winning edges $W$,
    assignment function $a: R\sqcup T \to \{0, \ldots, |T|\}$ 
\begin{algorithmic}[1]
    \State  $k \gets 1$; $A
    \gets \{\}$; $W \gets ()$; \Comment{Initialise} \label{line:Auc_Init} 
    \For{$r\in R$}
        \State $A
        (r) \gets \{\}$ \Comment{Initialise each robot's tasks}\label{line:Auc_InitAr} 
        \State $a(r)\gets 0$ \Comment{Assign each robot to itself} 
    \EndFor
    \For{$k\leq |T|$} \Comment{For each bid round}\label{line:Auc_LoopStart}
        \State{\edit{$w \gets \argmin c(s, t)$ \textbf{for} $s \in A\cup R$ \textbf{and} $t\in T\setminus A$}} \label{line:Auc_Calc}
        \State $W \gets W\cup (w)$ 
            \Comment{Append to list of winning edges}
            \label{line:Auc_RecordWinningEdge}
        \State $(w_0,w_1) \gets w$ \Comment{Label the endvertices of $w$}
        \State $a(w_1) \gets k$ \Comment{Label winning task with the round}\label{line:Auc_RecordWinningTask}
        \State $A
        \gets A
        \cup \{w_1\}$ \Comment{Update the assigned tasks}\label{line:Auc_AllAllocatedTasks}
        \For {$r\in R$} \Comment{Update the winning robot's tasks}
             \State \textbf{if} $w_0\in A
             (r)\cup \{r\}$ 
             \textbf{then} $A
             (r) \gets A
             (r) \cup \{w_1\}$ \label{line:Auc_RecordWinningRobot}
        \EndFor
     \State $k \gets k+1$ \Comment{Move onto the next bid round}
    \EndFor
    \State \Return $W, a$
\end{algorithmic}\label{alg:SSI_Auction}
\end{algorithm}

The \emph{sequential single-item (SSI)
auction} $\textsc{Auction}(\mathcal{M})$ on $\mathcal{M} = (R,T,L,p,c)$ which forms the first step of $\textsc{Assign}(\mathcal{M})$ consists of $|T|$ \emph{bid rounds} (see Algorithm \ref{alg:SSI_Auction}). 
At the start of the auction, 
the set of tasks $A(r)\subseteq T$ \emph{assigned to each robot} $r\in R$ is empty, and thus the set of \emph{all assigned tasks} $A=\bigsqcup_{r\in R} A(r) = \emptyset$ (see lines \ref{line:Auc_Init} and \ref{line:Auc_InitAr}).

In each \emph{bid round} $k\in\{1,\ldots, |T|\}$, every element in $A\sqcup R$ bids on every \emph{unassigned task} in $T\setminus A$. The pair $(s,t)\in (A\sqcup R) \times (T\setminus A)$ with the lowest \emph{bid} $c(s,t)$ wins the round (line \ref{line:Auc_Calc}). This pair is called the \emph{winning edge} $w$ as it corresponds to an edge $\{s,t\}$ in the robot-task graph $G(\mathcal{M})$. The task $t$ is called the \emph{winning task}. The \emph{winning robot} is the unique $r\in R$ for which $s\in A(r)$. We \emph{assign} $t$ to $r$ by appending it to $A(r)$ and $A$  (lines \ref{line:Auc_RecordWinningRobot} and  \ref{line:Auc_AllAllocatedTasks}), and record the round in which this assignment occurred by updating the \emph{assignment function} $a$  (line \ref{line:Auc_RecordWinningTask}). The winning edge is recorded as the $k^{\text{th}}$ element of the list $W$ (line \ref{line:Auc_RecordWinningEdge}).

The auction terminates when all tasks have been assigned ($A = T$), and outputs the list of winning edges $W$ and the assignment function $a$.

Note that there are numerous SSI-auctions that can be defined on $\mathcal{M}$ by changing the bid function. For an overview of several such SSI-auctions, 
see \cite{L+_AuctionBasedMRR}. 




\subsection{Step 2: \textsc{DFShortcut}} \label{subsec:DoubleTreeAlgo}

\begin{algorithm}[b]
\caption{\textsc{DFShortcut}}\label{alg:DepthFirstShortcut}
\textbf{input:} 
$R,T$, $\textsc{Auction}(\mathcal{M}) = (W,a)$\\
\textbf{output:} plan $\mathcal{P} = \{P(r): 
r\in R\}$ 
\begin{algorithmic}[1]
    \State $\mathcal{P} \gets \{\}$ \Comment{Initialise the plan}
    \For{$r\in R$} \Comment{Initialise tree and robot-route}
        \State $V(r) \gets \{r\}$; $P(r)\gets (r)$ 
    \EndFor
    \For{$w = (s_1,s_2)\in W$} 
        \For{$r\in R$} \Comment{Construct vertex set for each tree}
            \State \textbf{if} $s_1\in V(r)$ \textbf{then} $V(r) \gets V(r) \cup \{s_2\}$ 
        \EndFor
    \EndFor
    \For{$r\in R$} \Comment{Construct robot-route from tree}
        \While {$|P(r)| < |V(r)|$} 
            \State $a \gets \infty; v \gets \text{Null}$ \Comment{Initialise next task} 
            \For{$t \in P(r)$ \textbf{in reverse order}} 
                \For{$s\in V(r) \setminus \{u:u\in P(r)\}$}
                    \If{$(t,s)\in W \textbf{ or } (s,t)\in W$}
                        \If{$a(s) < a$}
                            \State $a \gets a(s)$ \Comment{Update minimum}
                            \State $v \gets s$ \Comment{Update next vertex}
                        \EndIf
                    \EndIf
                \EndFor
            \If{$v\neq\text{Null}$} \Comment{If next vertex known}
                \State $P(r) \gets P(r) \cup (v)$ \Comment{Append  vertex} 
                \State \textbf{break} \Comment{Return to the while loop}
            \EndIf \Comment{Otherwise, repeat the for loop} 
            \EndFor
        \EndWhile
        \State{$\mathcal{P}\gets \mathcal{P} \cup \{P(r)\}$} \Comment{Append robot-route to plan}
    \EndFor
    \State \Return $\mathcal{P}$
\end{algorithmic}    
\end{algorithm}


In the second step of \textsc{Assign}, we use \textsc{DFShortcut} (Algorithm \ref{alg:DepthFirstShortcut}) to convert the set of winning edges $W$ from $\textsc{Auction}(\mathcal{M})$ into robot-routes. Observe that 
each connected component $W_i$ of $W$ is the edge set of a tree $F_i=(V_i,W_i)$ in $G(\mathcal{M})$ which contains exactly one robot $r$ in its vertex set $V_i$.
\textsc{DFShortcut} uses this property to perform a depth-first search of each $F_i$, starting from $r$. The resulting vertex ordering gives a robot-route $P(r)$ starting at $r$ and visiting every other vertex in $V_i$ exactly once.

\subsection{2-approximation for metric \textsc{MinSum}} \label{subsec:2Approx} 

\begin{lem} \label{lem:DoubleTreeCost}
Let $\mathcal{M} = (R,T,L,p,c)$ be an MRTA instance with robot-task graph $G$, where $c$ is metric and injective on $E(G)$. Then 
$
\edit{C}(\textsc{Assign}(\mathcal{M})) \leq 2 \edit{C}(\textsc{MinSum}(\mathcal{M})).
$
\end{lem}
This is easy to see; we include the proof for completeness. 
To simplify notation, for any edge set $F$ 
we let $\edit{C}(F) = \sum_{f\in F} c(f)$ \edit{and $C'(F)= \sum_{f\in F} c'(f)$}. When 
$F$ is the edge set of a route $P$, this gives $\edit{C}(F)\leq \edit{C}(P)$ \edit{(and similarly $\edit{C}'(F)\leq \edit{C}'(P)$)}. If $P$ is a trail, this holds with equality.
\begin{pf}
Let 
$W$ denote the list of winning edges calculated by $\textsc{Auction}(\mathcal{M})$. 
We first convert $G$ to a complete graph $G'$ with edge-labelling $c'$ as follows. If for an edge $\{x,y\}\in E(G')$, we have $x,y\in R$ then $c'(x,y) = 0$, 
otherwise $c'(x,y) = c(x,y)$. 
Let $F'\subset E(G')$ be a spanning tree of $R$. Then both $W'=F'\cup W$ and $P'=F'\cup\{e\in E(P_i^*): P_i^*\in \textsc{MinSum}(\mathcal{M})\}$ are spanning trees of $G'$ with $\edit{C}'(W') = \edit{C}(W)$ and $\edit{C}'(P') = \edit{C}(\textsc{MinSum}(\mathcal{M}))$ respectively. 

Observe that $W'$ can be obtained by applying Prim's algorithm \cite{Prim} to $G'$. 
Thus $W'$ is a minimal spanning tree of $G'$, and so $\edit{C}'(W')\leq \edit{C}'(P')$. Combining this with the above two equalities gives
\begin{equation}\label{eqn:2approx_MST}
\edit{C}(W) \leq \edit{C}(\textsc{MinSum}(\mathcal{M})).
\end{equation}
The double-tree method $\textsc{DFShortcut}$ (or indeed any alternative double-tree short-cutting method) guarantees that for any connected component $W_i\subseteq W$, the corresponding robot-route $P_i\in\textsc{Assign}(\mathcal{M})$ satisfies $\edit{C}(P_i) \leq 2 \edit{C}(W_i)$. Thus
$
\edit{C}(\textsc{Assign}(\mathcal{M}))\leq 2\edit{C}(W).
$
Combining this with \eqref{eqn:2approx_MST}, we obtain
$
\edit{C}(\textsc{Assign}(\mathcal{M})) \leq 2\edit{C}(\textsc{MinSum}(\mathcal{M}))
$
as required.\hfill$\square$
\end{pf}

\section{Sensitivity analysis} \label{sec:Sensitivity}

Let $G(\mathcal{M})=(R,T;E;c)$ and $G(\mathcal{M}')=(R,T;E;c')$ be robot-task graphs for respective MRTA instances $\mathcal{M}$ and $\mathcal{M}'$, 
where $c'$ is obtained from $c$ by adjusting the cost of each edge $e$ 
of the robot-task graph $G(\mathcal{M})$. How large can these adjustments be if we wish to ensure that $\textsc{Assign}(\mathcal{M}')$ outputs the same plan as $\textsc{Assign}(\mathcal{M})$?  
By the definition of \textsc{DFShortcut}, 
this is guaranteed when $\textsc{Auction}(\mathcal{M})$ and $\textsc{Auction}(\mathcal{M}')$ output the same 
winning edges in the same order. 

Let $\delta(e)\in \mathbb{R}$ be such that $c'(e) = c(e)+\delta(e)$ for all $e\in E$. In this section, we 
obtain bounds $\underline{\Delta}(e),\overline{\Delta}(e)\in\mathbb{R}_+$ 
which guarantee $\textsc{Assign}(\mathcal{M}')$ outputs the same plan as  $\textsc{Assign}(\mathcal{M})$ so long as $\delta(e)\in(-\underline{\Delta}(e),\overline{\Delta}(e)]$ for all $e\in E$. In Subsection \ref{subsec:OneEdge}, we consider the case where $c'$ is obtained from $c$ by changing the cost of exactly one edge; and in Subsection \ref{subsec:AllEdges}, we allow the cost of all edges to change. Finally, in Subsection \ref{subsec:AllEdges_Best}, we show how to select the lexicographic best set of bounds from the candidates found in Subsection \ref{subsec:AllEdges}. First, we formalise the above ideas. 

\subsection{Families of intervals}\label{subsec:robust_notation}

Following the notation in Section \ref{sec:Auction}, we let $\textsc{Assign}(\mathcal{M}) = \mathcal{P}=\{P_i\}$ and $\textsc{Auction}(\mathcal{M}) = (W=(w_k),a)$. To avoid ambiguity, we append a superscript to the corresponding properties for $\mathcal{M}'$. So 
$\textsc{Assign}(\mathcal{M}') = \mathcal{P}'$ and $\textsc{Auction}(\mathcal{M'}) = (W'=(w'_k),a')$.

Recall that each bid round of $\textsc{Auction}(\mathcal{M})$ starts with a set of assigned tasks $A$, considers the set of all edges $e$ of $G(\mathcal{M})$ which have one end-point in $R\sqcup A$ and the other in $T\setminus A$, and selects the lowest cost edge in this set. For each edge $e=\{x,y\}$ we denote the set of bid rounds of $\textsc{Auction}(\mathcal{M})$ which consider $e$ by $B_e$. 
This 
can be determined from the assignment function $a$ by
\begin{equation}\label{eqn:Be_definition}
B_e = \left\{k\in \mathbb{Z}\ \middle\vert \begin{array}{l}
\min(a(x), a(y)) < k \text{ and} \\
k \leq \max(a(x),a(y))
\end{array}\right\}.
\end{equation}
%
For each bid round $k$, we denote the \emph{edge with the second-best bid} by $u_k$. In other words, $u_k$ is an edge with $k\in B_{u_k}$ and $u_k\neq w_k$, such that for all edges $e\edit{\neq w_k}$ with $k\in B_e$, we have $c(w_k)\leq c(u_k)\leq c(e)$. We denote the list of all second-best edges by $U=(u_k)$.

Our goal in this section is to find an interval in $\mathbb{R}$ for each edge $e$ of $G(\mathcal{M})$, such that for any fluctuation of the cost of $e$ within this interval, the outcome of the auction is maintained. More formally, for $G(\mathcal{M})=(R,T;E;c)$, 
\edit{let $\overline{\Delta}$ denote an \emph{upper bound} function and $\underline{\Delta}$ denote a \emph{lower bound} function, where
\begin{align}
  &\overline{\Delta}: E \to \mathbb{R}_+ \text{ and}  \nonumber \\
  &\underline{\Delta}: E \to \mathbb{R}_+ \text{ such that } \underline{\Delta}(e)\leq c(e) \text{ for all } e\in E.  \label{eqn:defn_RobustFamily}
\end{align}
Let $I$ be a function which maps each edge $e\in E$ to the \emph{interval} $(-\underline{\Delta}(e), \overline{\Delta}(e)]\subseteq \mathbb{R}$. The resulting image $I(E)$ is called}
a \emph{family of intervals}. We say $\textsc{Auction}(\mathcal{M})$ is \emph{robust} on $I\edit{(E)}$ \edit{(or equivalently on $I$)} 
if whenever 
$G(\mathcal{M}')=(R,T;E;c')$ satisfies $c'(e)\in (c(e)-\underline{\Delta}(e), c(e)+\overline{\Delta}(e)]$ for all $e\in E$, we have $w'_k = w_k$ for all bid rounds $k$. In Subsection \ref{subsec:AllEdges}, we show how to construct such a family $I\edit{(E)}$ for $\textsc{Auction}(\mathcal{M})$.
%
%

\subsubsection{The maximal family of intervals}

Given an MRTA instance $\mathcal{M}$ 
with robot-task graph $G(\mathcal{M})=(R,T;E;c)$, let $\mathcal{I}_\mathcal{M}$ represent the collection of all families $I$ of intervals such that $\textsc{Auction}(\mathcal{M})$ is robust on $I$.
We define $\rho : \mathcal{I}_\mathcal{M}\to \mathbb{R}_+^{2|E|}$ to be a 
function that given any $I=\{(-\underline{\Delta}(e),\overline{\Delta}(e)]: e\in E\}$, first concatenates all pairs $\{\underline{\Delta}(e),\overline{\Delta}(e)\}$ to form a vector $v\in \mathbb{R}_+^{2|E|}$ and then sorts $v$ into increasing order, giving $\rho(I)=v^\uparrow$. 

The \emph{lexicographic ordering} on $\mathbb{R}_+^{2|E|}$ is given by the relation $\succeq$, where for $u=(u_i),v=(v_i) \in \mathbb{R}_+^{2|E|}$, we say $u \succeq v$ if either $u_i=v_i$ for all $i\in\{1,\ldots, {2|E|}\}$, or there exists some $j\in\{1,\ldots, {2|E|}\}$ such that $u_j > v_j$ and $u_i=v_i$ for all $i<j$. If the latter case holds, we write $u \succ v$.

In Subsection \ref{subsec:AllEdges_Best}, we show how to find a family ${I}\in\mathcal{I}_\mathcal{M}$ such that $\rho({I}) \succeq \rho(I')$ for all $I'\in\mathcal{I}_\mathcal{M}$. We say such a $\rho({I})$ is \emph{maximal under the lexicographic ordering} of $\rho(\mathcal{I}_\mathcal{M})$, or equivalently that ${I}$ is \emph{maximal} in $\mathcal{I}_\mathcal{M}$. A family ${I}\in\mathcal{I}_\mathcal{M}$ is the \emph{unique maximal element} of $\mathcal{I}_\mathcal{M}$ if and only if $\rho({I}) \succ \rho(I')$ for all $I'\in\mathcal{I}_\mathcal{M}\setminus \{{I}\}$.

\subsection{Changing the cost of a single edge}\label{subsec:OneEdge}

We split our analysis into two cases: where the chosen edge wins a round of 
$\textsc{Auction}(\mathcal{M})$ (Lemma \ref{lem:ChangeOneEdge_winning}), and where it loses every round 
of $\textsc{Auction}(\mathcal{M})$ (Lemma \ref{lem:ChangeOneEdge_losing}).

\begin{lem}\label{lem:ChangeOneEdge_losing}
Let $G(\mathcal{M})=(R,T;E;c)$ and $G(\mathcal{M}')=(R,T;E;c')$ be two robot-task graphs where $c$ and $c'$ are injective on $E$. 
Let $\textsc{Auction}(\mathcal{M}) = (W=(w_k),a)$ and $\textsc{Auction}(\mathcal{M}') = (W'=(w'_k),a')$. Suppose there exists some 
$e\in E \setminus W$ such that $c'(f)=c(f)$ for all $f\in E\setminus \{e\}$.
Then $\textsc{Auction}(\mathcal{M}) = \textsc{Auction}(\mathcal{M}')$ if and only if 
$c'(e)>\max_{k\in B_{e}} c(w_{k})$.
\end{lem}

\begin{pf}
%
First, suppose $c'(e)>\max_{k\in B_{e}} c(w_{k})$. Then for each $k\in B_{e}$ we have
$c'(e) >c(w_{k}) = c'(w_{k})$. Thus $w'_{k} = w_{k}$ for all $k\in B_{e}$, 
which implies $w'_{k} = w_{k}$ for all rounds $k$. The definition of the assignment function then gives $a'=a$, and so $\textsc{Auction}(\mathcal{M}) = \textsc{Auction}(\mathcal{M}')$, as required.

Instead, suppose $c'(e)\leq \max_{k\in B_{e}} c(w_{k})$. Then there exists some minimal $K\in B_{e}$ such that $c'(e)\leq c(w_{K}) = c'(w_{K})$. If this holds with equality, then this violates the requirement that $c'$ is injective. 
Hence $c'(e) < c'(w_{K})$, 
which implies $w'_{K} =e \neq w_{K}$, and thus $\textsc{Auction}(\mathcal{M}) \neq \textsc{Auction}(\mathcal{M}')$.\hfill$\square$
\end{pf}

\begin{lem}\label{lem:ChangeOneEdge_winning}
Let $G(\mathcal{M})=(R,T;E;c)$ and $G(\mathcal{M}')=(R,T;E;c')$ be two robot-task graphs where $c$ and $c'$ are injective on $E$.
Let $\textsc{Auction}(\mathcal{M}) = (W=(w_k), a)$, $\textsc{Auction}(\mathcal{M}') = (W'=(w'_k), a')$, and $U=(u_k)$ denote the edges with second-best bids on $\textsc{Auction}(\mathcal{M})$.
Suppose there exists some 
$e=w_{K}\in W$ such that $c'(f) = c(f)$ for all $f\in E\setminus \{e\}$. 
Then $\textsc{Auction}(\mathcal{M}) = \textsc{Auction}(\mathcal{M}')$, if and only if
\[
c(u_{K})> c'(e)  \begin{cases}
> \max_{k\in B_{e}\setminus \{K\}} c(w_{k})
  &\hspace{0.2em}\text{if }B_{e}\setminus \{K\}\neq \emptyset, \\
\geq 0
  &\hspace{0.2em}\text{otherwise.}
\end{cases}
\]
\end{lem}
\begin{pf}
%
Observe that proving $\textsc{Auction}(\mathcal{M}) = \textsc{Auction}(\mathcal{M}')$ is equivalent to showing that $w'_k = w_k$ for all rounds $k$, since the assignment function is determined by the set of winning edges. We split the proof into the cases $c'(e)\leq c(e)$ and $c'(e)> c(e)$.

First suppose $c'(e)\leq c(e)$. Then $c(u_K) > c'(e)$, and it only remains to consider 
the lower bounds. If $B_{e}\setminus\{K\}\neq\emptyset$, then necessity and sufficiency follows from a similar argument to Lemma \ref{lem:ChangeOneEdge_losing}, with the small change that since $e=w_{K}$, we omit round $K$ when taking the maximum. 
So suppose $B_{e}=\{K\}$. Then $e$ is considered in exactly one round $K$ of $\textsc{Auction}(\mathcal{M})$, which it wins. Since $c'(e)\leq c(e)$, this implies $w'_K = e = w_K$, and thus $w'_{k} = w_{k}$ for all $k$. The lower bound $c'(e) \geq 0$ follows from the fact $c'$ is a cost function, and so is non-negative by definition.

Instead suppose $c'(e)>c(e)$. In this case, the lower bounds are trivial, and it only remains to consider the upper bound. To prove sufficiency, suppose $c(u_K)>c'(e)$. 
Since $c'(e)> c(e)$, 
we have $w'_{k} = w_{k}$ for all $k\leq K$, 
which implies $w'_{k}=w_{k}$ for all rounds $k$.
To show necessity, suppose $c'(e)\geq c(u_K)$. This cannot hold with equality since 
$c'$ is injective. So we must have $c'(e) > c(u_K) = c'(u_K)$, which implies $u_K=w'_{K}\neq w_{K}$.\hfill$\square$
\end{pf}
Lemmas \ref{lem:ChangeOneEdge_losing} and \ref{lem:ChangeOneEdge_winning} immediately imply the following.

\begin{cor}\label{cor:ChangeOneEdge_Summary}
Let $\mathcal{M}=(R,T,L,p,c)$ and $\mathcal{M'}=(R,T,L,p,c')$ be MRTA instances whose cost functions $c$ and $c'$ are injective on the edge sets of the respective robot-task graphs $G(\mathcal{M})$,  and $G(\mathcal{M}')$, 
and differ only in the cost of a single edge $e$. Let $\textsc{Auction}(\mathcal{M}) = (W=(w_k), a)$, and let $U=(u_k)$ denote the edges with second-best bids on $\textsc{Auction}(\mathcal{M})$. 
Then $\textsc{Assign}(\mathcal{M}') = \textsc{Assign}(\mathcal{M})$ if either
\begin{enumerate}
\item $e\not\in W$ and $c'(e)> \max_{k\in B_{e}} c(w_{k})$,
\item $e =w_{K}\in W$ with $|B_{e}|=1$, and $c(u_{K}) > c'(e) \geq 0$, or \label{part:OneEdge_WinUniqueRound}
\item $e =w_{K}\in W$ with $|B_{e}|>1$, and $c(u_{K}) > c'(e) > \max_{k\in B_{e}\setminus\{K\}} c(w_{k})$.
\end{enumerate}
\end{cor}


\subsection{Changing the cost of all edges}\label{subsec:AllEdges}


Let $G(\mathcal{M}) = (R,T; E; c)$ with $c$ injective on $E$, and $\textsc{Auction}(\mathcal{M}) = (W=(w_k), a)$ with $U=(u_k)$ denoting the list of edges with second-best bids on $\textsc{Auction}(\mathcal{M})$.  
Lemma \ref{lem:ChangeOneEdge_winning} implies that if $\textsc{Auction}(\mathcal{M})$ is robust on a family  of intervals  $\{(-\underline{\Delta}(e),\overline{\Delta}(e)] : e \in E\}$, then 
\begin{equation}\label{eqn:OneEdge_WinCondition}
0 \leq \overline{\Delta}(w_k) < c(u_k)-c(w_k) \text{ for all } k\in\{1,\ldots, |T|\}.
\end{equation}
In this subsection, we build on this observation by showing that any set $I_0(W) = \{\overline{\Delta}(w_k) : w_k \in W\}$ which satisfies \eqref{eqn:OneEdge_WinCondition} can be extended to an interval family $I(E)=\{(-\underline{\Delta}(e),\overline{\Delta}(e)] : e \in E\}$ using the \textsc{ErrorIntervals} algorithm (see below). Theorem \ref{thm:RobustMargins_general} shows that $\textsc{Auction}(\mathcal{M})$ is robust on the resulting family 
$I(E)$. Observe that lines \ref{line:robust_win}, \ref{line:robust_win_UniqueBidRound} and \ref{line:robust_lose} of \textsc{ErrorIntervals} are inspired by the inequalities in Lemmas  \ref{lem:ChangeOneEdge_losing} and \ref{lem:ChangeOneEdge_winning}, and  line \ref{line:robust_BidSets} comes from equation \eqref{eqn:Be_definition}.

Theorem \ref{thm:RobustMargins_BestExtension} then shows that \textsc{ErrorIntervals} calculates 
the best possible family of intervals containing $I_0(W)$ as upper bounds 
on the winning edges. Here ``best possible'' means that if we expand any interval, then $\textsc{Auction}(\mathcal{M})$ is not robust on the resulting interval family. This result will be key to obtaining the maximal interval family for $\textsc{Auction}(\mathcal{M})$ in Subsection \ref{subsec:AllEdges_Best}. 

\begin{algorithm}
    \caption{\textsc{ErrorIntervals}
    }\label{alg:RobustMargins}
    \textbf{input:} 
    $G=(R,T;E;c)$, 
    $W=(w_k)$, 
    $a: R \sqcup T \to \{0, \ldots, |T|\}$,  
    $I_0 : W\to \mathbb{R}_+$  
    \\
    \textbf{output:} 
    function $I(e) =(-\underline{\Delta}(e), \overline{\Delta}(e)]$ for all $e\in E$, where $\underline{\Delta}(e),\overline{\Delta}(e)\in \mathbb{R}_+$ 
    
    \begin{algorithmic}[1]
    \For{$w_k\in W $}
        \State $\overline{\Delta}(w_k) \gets I_0(w_k)$ \Comment{Initialise}
       \label{line:robust_win_I0}
    \EndFor
    \For{$e = \{x,y\}\in E$}
        \State $\mu_1\gets\min(a(x), a(y))$; $\mu_2\gets\max(a(x), a(y))$
        \State $B_e\gets \{k\in \mathbb{Z}: \mu_1 < k \leq \mu_2\}$ \label{line:robust_BidSets}
        \If{$\exists K$ with $e = w_K\in W$ and $B_e\setminus\{K\}\neq\emptyset$}
            \State $\underline{\Delta}(e)\gets c(e) - \max_{k\in B_e \setminus \{K\}} c(w_k) + \overline{\Delta}(w_k)$ \label{line:robust_win}
        \ElsIf{$\exists K$ with $e = w_K\in W$}
            \State $\underline{\Delta}(e)\gets c(e)$ \label{line:robust_win_UniqueBidRound}
        \Else
            \State $\overline{\Delta}(e) \gets \infty$ \label{line:robust_lose_infty}
            \State $\underline{\Delta}(e) \gets c(e) - \max_{k\in B_e} c(w_k) + \overline{\Delta}(w_k)$ \label{line:robust_lose}
        \EndIf
        \State $I(e) \gets (-\underline{\Delta}(e), \overline{\Delta}(e)]$
    \EndFor
    \State \Return $I$
\end{algorithmic}
\end{algorithm}

\begin{thm}\label{thm:RobustMargins_general}
Let $G(\mathcal{M})=(R,T;E;c)$ be a robot-task graph with $c$ injective on $E$, $\textsc{Auction}(\mathcal{M}) = (W=(w_k), a)$, and $U=(u_k)$ be the list of edges with second-best bids on $\textsc{Auction}(\mathcal{M})$. Let $I_0: W \to \mathbb{R}_+$ such that 
$0\leq I_0(w_k) <c(u_k) - c(w_k)$ for all 
$k$. 
If $I = \textsc{ErrorIntervals}(G(\mathcal{M}), W, a, I_0)$, then $\textsc{Auction}(\mathcal{M})$ is robust on the family 
$
I(E).
$
\end{thm}
\begin{pf}
Denote $I(E) = \{(- \underline{\Delta}(e), \overline{\Delta}(e)]: e\in E\}$. 
Suppose $\mathcal{M}'=(R,T,L,p,c')$ is an MRTA instance with $c'(e)\in (c(e) - \underline{\Delta}(e), c(e)+\overline{\Delta}(e)]$ for all $e\in E$. Let $W'=(w'_k)$ denote the list of winning edges of $\textsc{Auction}(\mathcal{M}')$. It suffices to show that $w_k' = w_k$ for all $k$. 
We proceed by strong induction.

First suppose $k=1$. Since no tasks have yet been assigned, the first rounds of $\textsc{Auction}(\mathcal{M})$ and $\textsc{Auction}(\mathcal{M}')$ consider exactly the same set of edges; namely  $F_1=\{\{x,y\}\in E : x\in R \text{ and } y\in T\}$.  
Suppose $e\in F_1\setminus\{w_1\}$. Then $1\in B_e$.  So by the choice of $c'$, and line \ref{line:robust_lose} of \textsc{ErrorIntervals}, we have
\[
c'(e) > c(e) - \underline{\Delta}(e) 
\geq c(w_1) + \overline{\Delta}(w_1) \geq c'(w_1).
\]
Thus $w'_1=w_1$.

Now consider bid round $K>1$, and suppose that $w'_k = w_k$ for all $k<K$. Then the set of tasks that remain to be assigned in round $K$ is identical in both auctions. Thus exactly the same set of edges $F_K$ is considered  in round $K$ of $\textsc{Auction}(\mathcal{M})$ and $\textsc{Auction}(\mathcal{M}')$. Since $e\in F_K$ if and only if $K\in B_e$, a similar argument to the base case gives $w'_K = w_K$.\hfill$\square$
\end{pf}
The following result shows that $\textsc{ErrorIntervals}(\mathcal{M})$ calculates the 
maximal family corresponding to a given initialising set $I_0$.

\begin{thm}\label{thm:RobustMargins_BestExtension}
Let $G(\mathcal{M})=(R,T;E;c)$ be a robot-task graph with $c$ injective on $E$, and let $\textsc{Auction}(\mathcal{M}) = (W=(w_k),a)$ with edges $ U=(u_k)$ of second-best bids. Let $I_0: W \to \mathbb{R}_+$ with $0\leq I_0(w_k) < c(u_k)-c(w_k)$ for all $k$, and $I= \textsc{ErrorIntervals}(G(\mathcal{M}),W,a, I_0)$. Let $\mathcal{I}_0$ be the collection of all interval families $I'(E) = \{(- \underline{\Delta}'(e), \overline{\Delta}'(e)]: e\in E\}$ which satisfy
$\overline{\Delta}'(w_k)=I_0(w_k)$ for all $w_k\in W$, and such that $\textsc{Auction}(\mathcal{M})$ is robust on $I'(E)$. Then 
$I(E)$ is the unique maximal element of $\mathcal{I}_0$.
\end{thm}

\begin{pf}
It suffices to show $I'(e)\subseteq I(e)$ for all $e\in E$ and $I'\in\mathcal{I}_0$. We shall do this by letting  $I(E) = \{(- \underline{\Delta}(e), \overline{\Delta}(e)]: e\in E\}$, and proving that for all $I'(E) = \{(- \underline{\Delta}'(e), \overline{\Delta}'(e)]: e\in E\}\in\mathcal{I}_0$, we have both $\underline{\Delta}'(e) \leq \underline{\Delta}(e)$ and  $\overline{\Delta}'(e)\leq \overline{\Delta}(e)$ for all $e\in E$.

First consider the upper bounds. 
By the theorem hypotheses and line \ref{line:robust_lose_infty} of \textsc{ErrorIntervals} respectively, we have $\overline{\Delta}'(e) = \overline{\Delta}(e)$ for all $e\in W$ and $\overline{\Delta}'(e) \leq \overline{\Delta}(e) = \infty$ for all $e\in E\setminus W$, as required.

We now prove the inequality for the lower bounds by contradiction. Assume that for some $e\in E$ we have $\underline{\Delta}'(e) > \underline{\Delta}(e)$. 
Pick some $\delta\in\mathbb{R}$ with $\underline{\Delta}'(e) > \delta > \underline{\Delta}(e)$, 
and let $G(\mathcal{M}')=(R,T;E,c')$  where $c'$ is the edge-labelling  given by $c'(e) = c(e) - \delta$ and $c'(f) = c(f) + \overline{\Delta}'(f)$ for all $f\in E\setminus \{e\}$. Then $
c(f) - \underline{\Delta}'(f) < c'(f) \leq c(f) + \overline{\Delta}'(f)
$ 
for all $f\in E$. Since $\textsc{Auction}(\mathcal{M})$ is robust on $I'(E)$, this implies $W' = W$ where $W'$ is the list of winning edges of $\textsc{Auction}(\mathcal{M}')$.

Suppose $e= w_K\in W$ and $B_e=\{K\}$. Then line \ref{line:robust_win_UniqueBidRound} of \textsc{ErrorIntervals} gives $\underline{\Delta}(e) = c(e)$. So the assumption $\underline{\Delta}'(e) > \underline{\Delta}(e)$ implies $\underline{\Delta}'(e) > c(e)$, which contradicts 
\eqref{eqn:defn_RobustFamily}. 
Instead suppose either $e\in E\setminus W$ or $|B_e|\geq 2$. Then by lines \ref{line:robust_lose} or \ref{line:robust_win} respectively of \textsc{ErrorIntervals}, there is some minimal $k\in B_e$ such that
\[
c(e) - \underline{\Delta}(e) = 
c(w_k) + \overline{\Delta}(w_k).
\]
Since $\delta > \underline{\Delta}(e)$ and $\overline{\Delta}'(w_k) = \overline{\Delta}(w_k)$, this implies
\[
c'(e) = c(e) - \delta < c(e) - \underline{\Delta}(e) = c'(w_k).
\]
Hence $e = w'_k\neq w_k$, which contradicts $W'=W$.\hfill$\square$
\end{pf}
There is still room for improvement. \textsc{ErrorIntervals} takes an initial set $I_0(W)$ of upper bounds on the winning edges of $\textsc{Auction}(\mathcal{M})$ as input. What  choice of $I_0(W)$ results in the maximal interval family? 
We answer this in the following subsection.

\subsection{The maximal interval family} \label{subsec:AllEdges_Best}

Let $\mathcal{M}$ be an MRTA instance with corresponding robot-task graph $G$. 
Observe that if $\textsc{Auction}(\mathcal{M})=(W,a)$ is robust on interval family $I(E)= \{(-\underline{\Delta}(e),\overline{\Delta}(e)]: e\in E(G)\}$ then this implies two properties of $I(E)$. Firstly, that every edge $e\in E(G)$ satisfies
\begin{equation}\label{eqn:LexBest_UnderbidderCondition_v2}
c(e) - \underline{\Delta}(e)\geq \max_{j\in B} c(w_j) + \overline{\Delta}(w_j),
\end{equation}
where $B = B_e\setminus\{K\}$ if $e=w_K\in W$, and $B=B_e$ otherwise. 
Rearranging this, we obtain $\overline{\Delta}(w_j) + \underline{\Delta}(e)\leq c(e)-c(w_j)$ for all edges $e\in E(G)\setminus\{w_j\}$ with $j\in B_e$. Since there are potentially many such edges $e$ for each $w_j$, this suggests that if $I(E)$ is the maximal element of $\mathcal{I}_\mathcal{M}$, then it \edit{is likely (but not guaranteed) that it also satisfies a} 
second property:
\begin{equation}\label{eqn:LexBest_FirstGuess_v2}
\overline{\Delta}(w_j) \leq \frac{1}{2} (c(e)-c(w_j))
\end{equation}
for all $e\in E(G)\setminus\{w_j\}$ with $j\in B_e$. 

In this section, we introduce the algorithm \textsc{Initialiser}, which uses these two observations to construct an initialising set \edit{$I_0(W)$} 
for $\textsc{ErrorIntervals}$. 
\edit{It does this by iteratively improving two estimates. For each losing edge $e$, we keep track of the largest lower bound $L(e)\leq c(e) - \underline{\Delta}(e)$ for inequality \eqref{eqn:LexBest_UnderbidderCondition_v2} found so far (line \ref{line:I0_UpdateConstraints}). And for each winning edge $w_k$, we keep track of the smallest upper bound $M(e) \geq c(w_k) + I_0(w_k) = c(w_k) + \overline{\Delta}(w_k)$ found so far over all edges $e$ that lose round $k$ (line \ref{line:I0_SetBound}). Inspired by inequality \eqref{eqn:LexBest_FirstGuess_v2}, we set $M(e)$ as the midpoint $\frac{1}{2}(c(w_k) +c(e))$ unless $L(e)$ already gives us a bigger upper bound (line \ref{line:I0_BestBound}).
}

Theorem \ref{thm:LexBestFamily} shows that \textsc{Initialiser} does indeed construct the initialising set 
for the unique maximal element of 
$\mathcal{I}_\mathcal{M}$.

\begin{algorithm}[h]
    \caption{\textsc{Initialiser}.}\label{alg:BestInitializer}
    \textbf{input:} 
    $G=(R,T;E;c)$,
    $W=(w_k)$,
    $a: R\sqcup T \to \{0,\ldots, |T|\}$
    \\
    \textbf{output:} $I_0 : W\to \mathbb{R}_+$
    \begin{algorithmic}[1]
    \For{$e\in E$}
        \State $L(e) \gets 0$ \Comment{Initialise lower bound in equation \eqref{eqn:LexBest_UnderbidderCondition_v2}}\label{line:I0_L0}
    \EndFor
    \State $S \gets$ \textbf{sort} $\{1, \ldots, |T|\}$ \textbf{by} $c(w_k)$ \textbf{decreasing} 
    \For{$k\in S$}  \Comment{Loop through $k$ in decreasing $c(w_k)$}
        \State $E_k\gets \{\}$ \Comment{Initialise}
        \For{ $e=\{x,y\}\in E$}
            \State $\mu_1 \gets \min(a(x),a(y))$ ; $\mu_2 \gets \max(a(x),a(y))$
            \If{ $\mu_1<k\leq\mu_2$}
                \If{ $(x,y)\not\in W$ \textbf{and} $(y,x)\not\in W $}
                    \State $E_k\gets E_k \cup e$ \Comment{Round $k$ considers $e$}
                \EndIf
            \EndIf
        \EndFor
        \For{$e\in E_k$} \Comment{For edges considered in round $k$}
            \State $M(e) \gets \max\left(L(e), \frac{c(w_k) + c(e)}{2}\right)$ 
            \label{line:I0_BestBound}
        \EndFor
        \State $I_0(w_k) \gets \min_{e\in E_k} M(e) - c(w_k)$ \label{line:I0_SetBound}
        \For{$e\in E_k$}
            \State $L(e) \gets \max(L(e), c(w_k) + I_0(w_k))$ 
            \label{line:I0_UpdateConstraints}
        \EndFor
    \EndFor
    \State \Return $I_0$
    \end{algorithmic}
\end{algorithm}

\begin{thm}\label{thm:LexBestFamily}
Let $\mathcal{M}$ be an MRTA instance with robot-task graph  
$G = (R,T;E;c)$ such that $c$ is injective on $E$. Let $\textsc{Auction}(\mathcal{M}) = (W=(w_k),a)$, and $I_0 = \textsc{Initialiser}(G,W,a)$. Further, let $I = \textsc{ErrorIntervals}(G,W,a,I_0)$.
Then 
$I(E)$ is the unique maximal element of $\mathcal{I}_\mathcal{M}$.
\end{thm}

\begin{pf}
Assume for a contradiction that there exists some $I'(E)\in\mathcal{I}_\mathcal{M}\setminus \{I(E)\}$ with $\rho(I'(E))\succeq \rho(I(E))$. We denote $I(E)=\{(-\underline{\Delta}(e),\overline{\Delta}(e)]: e\in E\}$ and $\rho(I(E)) = (\rho_i)$. 
Similarly, we write $I'(E)=\{(-\underline{\Delta}'(e),\overline{\Delta}'(e)]: e\in E\}$ and $\rho(I'(E)) = (\rho'_i)$.
%

Since $I\neq I'$, 
there is some $e\in E$ and symbol $\Delta\in\{\underline{\Delta}, \overline{\Delta}\}$ such that $\Delta(e)\neq \Delta'(e)$ (so either $\underline{\Delta}(e)\neq \underline{\Delta}'(e)$ or $\overline{\Delta}(e)\neq \overline{\Delta}'(e)$). 
Pick such an $e$ with $\Delta(e)\in\mathbb{R}_+$ minimal. 
By relabelling identical values in $\rho(I(E))$, we may assume there exists some $j$ such that $\rho_j = \Delta(e)$, and
$\rho_i = \rho'_i$ for all $i<j$. 
Since $\rho(I'(E))\succeq \rho(I(E))$ and $\Delta(e)$ is minimal, this implies 
\begin{equation}\label{eqn:LexBestPf_BoundInequality}
\rho_j = \Delta(e) < \Delta'(e) = \rho'_{j'} \quad\text{ for some } j'\geq j;
\end{equation}
and for all $i< j$, there exists some $e_i\in E$ such that either
\begin{align}
 \rho_i = \overline{\Delta}(e_i) &= \overline{\Delta}'(e_i) = \rho'_i \quad\text{or} \label{eqn:LexBestPf_RhoEntries_upper}\\
 \rho_i = \underline{\Delta}(e_i) &= \underline{\Delta}'(e_i) =  \rho'_i. \label{eqn:LexBestPf_RhoEntries_lower}
\end{align}
Our goal is to show that such a $\Delta(e)$ does not exist. 
The proof splits into cases depending on whether $\Delta$ corresponds to $\overline{\Delta}$ or $\underline{\Delta}$. We first prove a useful fact.

\begin{claim}\label{Claim:UpperBound_f}
If $\Delta(e)$ corresponds to $\overline{\Delta}(e)$, then $e=w_K$ for some $w_K\in W$, and there exists some edge $f\neq w_K$ with $K\in B_f$ such that $c(w_K) + \overline{\Delta}(w_K) = c(f) - \underline{\Delta}(f)$.
\end{claim}

\begin{pf}
By \eqref{eqn:LexBestPf_BoundInequality} and the definition of \textsc{ErrorIntervals}, we have $\overline{\Delta}(e)< \infty$ 
which implies $e = w_K\in W$ for some $K$, 
giving the first part of the claim. 
In \textsc{Initialiser}, let $L_0$ denote the value of the function $L$ in line \ref{line:I0_L0}, and $L_k$ (respectively $M_k$) denote the value of $L$ (resp.\ $M$) given in line \ref{line:I0_UpdateConstraints} (resp. line \ref{line:I0_BestBound}) for bid round $k\in S$. 
By the definition of \textsc{Initialiser}, there exists some $f\in E_K$ such that
\begin{equation}\label{eqn:LexBestPf_Claim1_wk}
c(w_K) + I_0(w_K) = M_K(f) = L_K(f);
\end{equation}
where the first equality follows from line \ref{line:I0_SetBound}, and the second from lines \ref{line:I0_BestBound} and \ref{line:I0_UpdateConstraints}. We shall show that this $f$ is the edge sought.

Suppose $L_K(f)$ is used as input in line \ref{line:I0_BestBound} for some subsequent $k\in S$. 
Then, by the ordering of $S$, we have $c(w_k)\leq c(w_K)$. Hence 
\begin{align}\label{eqn:LexBestPf_Claim1_Mf}
M_k(f) &= \max\left(L_K(f), \frac{c(w_k) + c(f)}{2}\right)\\
       &\leq \max\left(L_K(f), \frac{c(w_K) + c(f)}{2}\right) = L_K(f). \notag
\end{align}
Thus, applying lines \ref{line:I0_UpdateConstraints} and \ref{line:I0_SetBound} of \textsc{Initialiser} and equation \eqref{eqn:LexBestPf_Claim1_Mf} gives
\begin{align*}
L_k(f) &= \max(L_K(f), c(w_k)+I_0(w_k)) \\
       &= \max\left(L_K(f), \min_{e'\in E_k} M(e'))\right)\\
       &\leq \max(L_K(f), M_k(f)) = L_K(f). 
\end{align*}
Hence $L_k(f) = L_K(f)$, and so the value of $L(f)$ remains unchanged for any subsequent iterations over $S$. Thus, when \textsc{Initialiser} completes, we have 
\begin{equation}\label{eqn:LexBestPf_Claim1_Lf} 
L_K(f)= L(f) = \max_{i\in B} c(w_i) + I_0(w_i), 
\end{equation}
where $B=B_f\setminus\{j\}$ if $f=w_j$ for some $w_j\in W$, and $B=B_f$ otherwise. By line \ref{line:robust_win_I0} of \textsc{ErrorIntervals}, we have $\overline{\Delta}(w_i) = I_0(w_i)$ for all $i\in B$. 
Substituting this into line \ref{line:robust_win} of \textsc{ErrorIntervals} when $f\in W$, or line \ref{line:robust_lose} otherwise, gives
\begin{align}
c(f) - \underline{\Delta}(f) &= \max_{i\in B} c(w_i) + \overline{\Delta}(w_i) \label{eqn:LexBestPf_Claim1_cf} \\
 &= \max_{i\in B} c(w_i) + I_0(w_i).  \notag
\end{align}
Combining equations \eqref{eqn:LexBestPf_Claim1_wk}, \eqref{eqn:LexBestPf_Claim1_Lf}, and \eqref{eqn:LexBestPf_Claim1_cf} and using the fact $I_0(w_K) = \overline{\Delta}(w_K)$, we obtain
$
c(w_K) + \overline{\Delta}(w_K) = c(f) - \underline{\Delta}(f),
$
as required.\hfill$\square$
\end{pf}
We now proceed with the case analysis.

\smallskip\noindent 
\emph{Case 1.} $\Delta(e)$ corresponds to $\overline{\Delta}(e)$.\\
\noindent
By Claim \ref{Claim:UpperBound_f}, $e=w_K\in W$, and there exists an edge $f$ with $c(w_K) + \overline{\Delta}(w_K) = c(f) - \underline{\Delta}(f)$. 

First, suppose  $\underline{\Delta}(f) > \overline{\Delta}(w_K)$. 
By equation \eqref{eqn:LexBestPf_Claim1_wk} in the argument for Claim \ref{Claim:UpperBound_f}, and line \ref{line:I0_BestBound} of \textsc{Initialiser}, we know that $f$ satisfies
$
c(w_K) + I_0(w_K) = M_K(f) \geq \frac{1}{2}(c(w_K) + c(f)).
$ 
Since $I_0(w_K) = \overline{\Delta}(w_K)$ and $\underline{\Delta}(f) > \overline{\Delta}(w_K)$, this implies
$
c(w_K) + \overline{\Delta}(w_K) \geq c(f) - \overline{\Delta}(w_K) > c(f) - \underline{\Delta}(f),
$
which contradicts Claim \ref{Claim:UpperBound_f}.

So instead suppose $\underline{\Delta}(f) \leq \overline{\Delta}(w_K)$. 
Using the fact $\textsc{Auction}(\mathcal{M})$ is robust on $I'(E)$, together with equation \eqref{eqn:LexBestPf_BoundInequality} and  Claim \ref{Claim:UpperBound_f} we obtain
\begin{align*}
c(f) - \underline{\Delta}'(f) 
  &\geq c(w_K) +\overline{\Delta}'(w_K) 
  > c(w_K) + \overline{\Delta}(w_K) \\ 
  &= c(f) - \underline{\Delta}(f). 
\end{align*}
Thus $\underline{\Delta}'(f) < \underline{\Delta} (f)$. 
If $\underline{\Delta} (f) < \overline{\Delta} (w_K)$, then the fact $\underline{\Delta}'(f) \neq \underline{\Delta} (f)$ contradicts the choice of $\Delta(e)$.
So we must have $\underline{\Delta} (f) = \overline{\Delta} (w_K)$. In which case, equations \eqref{eqn:LexBestPf_BoundInequality}, \eqref{eqn:LexBestPf_RhoEntries_upper} and \eqref{eqn:LexBestPf_RhoEntries_lower} give 
$\rho_j = \underline{\Delta}(f) > \underline{\Delta}'(f) \geq \rho'_j$, which implies $\rho(I(E))\succ \rho(I'(E))$, contradicting our initial assumption.

\medskip \noindent 
\emph{Case 2.} $\Delta(e)$ corresponds to $\underline{\Delta}(e)$.\\
\noindent
Since $I'(E)\in\mathcal{I}_{\mathcal{M}}$, we have $c(e)-\underline{\Delta}'(e)\geq 0$ by 
\eqref{eqn:defn_RobustFamily}. This and equation \eqref{eqn:LexBestPf_BoundInequality} give $c(e) > \underline{\Delta}(e)$. 
By lines \ref{line:robust_win}, \ref{line:robust_win_UniqueBidRound} and \ref{line:robust_lose}  of \textsc{ErrorIntervals}, there is some $k\in B_e$ such that $e\neq w_k$ and
\begin{equation}\label{eqn:Lex_Lower_main}
    c(e) - \underline{\Delta}(e) = \max_{i\in B} c(w_i) + \overline{\Delta}(w_i) 
    = c(w_k) + \overline{\Delta}(w_k); 
\end{equation}
where $B = B_e\setminus\{K\}$ if $e=w_K\in W$, and $B=B_e$ otherwise. If there are multiple choices of round $k\in B$ which satisfy \eqref{eqn:Lex_Lower_main}, we choose one with $c(w_k)$ maximal.

First suppose $\overline{\Delta}(w_k)\leq \underline{\Delta}(e)$. Equation \eqref{eqn:LexBestPf_BoundInequality} and the fact  $\textsc{Auction}(\mathcal{M})$ is robust on  $I'(E)$ give
$
c(w_k)+\overline{\Delta}'(w_k)\leq c(e) - \underline{\Delta}'(e) < c(e) - \underline{\Delta}(e).
$
This and \eqref{eqn:Lex_Lower_main} imply $\overline{\Delta}'(w_k)< \overline{\Delta}(w_k)$. Combining this with conditions \eqref{eqn:LexBestPf_BoundInequality}, \eqref{eqn:LexBestPf_RhoEntries_upper} and \eqref{eqn:LexBestPf_RhoEntries_lower} gives 
$
\rho'_j \leq \overline{\Delta}'(w_k) 
< \overline{\Delta}(w_k)
\leq \underline{\Delta}(e) = \rho_j.
$
Thus $\rho(I(E))\succ \rho(I'(E))$, contradicting our initial assumption.

Instead suppose $\overline{\Delta}(w_k) > \underline{\Delta}(e)$. Line \ref{line:I0_SetBound} of \textsc{Initialiser} and line \ref{line:robust_win_I0} of \textsc{ErrorIntervals} imply 
\begin{equation}\label{eqn:MkeBound}
c(w_k) +\overline{\Delta}(w_k) \leq M_k(e).
\end{equation}
If $M_k(e) = \frac{1}{2}(c(w_k) + c(e))$, then substituting this into \eqref{eqn:MkeBound} gives 
$
c(w_k) +\overline{\Delta}(w_k)
< c(e) - \underline{\Delta}(e) 
$, which contradicts \eqref{eqn:Lex_Lower_main}.
Thus, by lines \ref{line:I0_BestBound} and \ref{line:I0_UpdateConstraints} of \textsc{Initialiser}, there exists some $l\in B$ with $c(w_l)> c(w_k)$ 
such that $M_k(e) = L_l(e) = c(w_l)+I_0(w_l)$. Using this, \eqref{eqn:Lex_Lower_main}, \eqref{eqn:MkeBound}, and the fact $I_0(w_l) = \overline{\Delta}(w_l)$, we obtain $c(w_l) + \overline{\Delta}(w_l) = c(e) - \underline{\Delta}(e)$, which contradicts the choice of $k$.\hfill$\square$ 
\end{pf}

The theory developed in this paper culminates in the following statement, which combines Lemma \ref{lem:DoubleTreeCost} with Theorem \ref{thm:LexBestFamily}. 

\begin{cor}
Let $\mathcal{M}$ be an MRTA instance with robot-task graph $G=(R,T;E;c)$ such that $c$ is metric and injective on $E$. Then
\begin{enumerate}
    \item $\textsc{Assign}(\mathcal{M})$ provides a 2-approximation for $\textsc{MinSum}(\mathcal{M})$, and
    \item $\textsc{Assign}(\mathcal{M})$ is robust over the interval family output by $\textsc{AuctionSensitivity}(\mathcal{M},W,a)$, where $(W,a) = \textsc{Auction}(\mathcal{M})$.
\end{enumerate}

\end{cor}  
\section{Simulations} \label{sec:Simulations}


The complexity of the algorithms can be readily derived and are summarised in Table \ref{tab:complexity}. These are seen to be relatively efficient in practice, and scale in a manner to allow implementation with even large numbers of agents and tasks.

\begin{table}[h]
    \caption{Complexity of  $\textsc{Assign}$, $\textsc{AuctionSensitivity}$ and their subroutines on a suitable MRTA instance $\mathcal{M} = (R,T,L,p,c)$. To simplify the statements, $m = |R|$ and $n= |T|$.\label{tab:complexity}}
    \centering
    \begin{tabular}{|l|c|}
    \hline
        Algorithm & Complexity \\
    \hline
         $\textsc{Auction}$ & $o(n^2(n+m))$ \\
         $\textsc{DFShortcut}$ & $o(m n^3)$\\
         \rowcolor{black!10} $\textsc{Assign}$ & $o(m n^3)$ \\ 
    \hline
         $\textsc{Initialiser}$ & $o(n(m+n)^2)$ \\
         $\textsc{ErrorIntervals}$ & $o(n(m+n)^2)$ \\
         \rowcolor{black!10}  $\textsc{AuctionSensitivity}$ & $o(n(m+n)^2)$ \\
    \hline  
    \end{tabular}        
\end{table}

\edit{These algorithms are efficient to run in practice. This was verified in 1000 simulations run in Matlab R2021a on a Intel i7 GPU running at 2.8 GHz where in each simulation an MRTA instance on 10 robots and 100 tasks with randomised positions and Euclidean distance as cost metric were considered. The average computation time for the subroutines was the following: 6.0ms for \textsc{Auction}, 0.5ms for \textsc{DFShortcut}, 737.7ms for \textsc{Initialise}, and 21.7ms for \textsc{ErrorIntervals}.}

We close this section with an example to illustrate situations where our sensitivity analysis would be helpful.

\begin{exmp}
Let $\mathcal{M} = (R,T,L,p,c)$ be the MRTA problem shown in Figure \ref{subfig:Example_InitialProb}, where $R=\{r_1, r_2\}$, $T=\{t_1,t_2, t_3\}$, and $p$ is given by the coordinates shown.
The grey box corresponds to an obstacle $B$ which must be avoided, and so the 
landscape we may traverse is $L = \mathbb{R}^2\setminus B$. For simplicity, we assume that for all $p_1,p_2\in L$, the cost $c(p_1,p_2)$ is given by the length of the shortest path between $p_1$ and $p_2$ that is contained entirely in $L$. 

Figure \ref{subfig:Example_InitialSoln} shows the plan output by $\textsc{Assign}(\mathcal{M})$. The upper and lower bounds for the interval family calculated by $\textsc{AuctionSensitivity}(\mathcal{M}, \textsc{Auction}(\mathcal{M}))$ are shown in Table \ref{table:example_sensitivity}.

\begin{figure}[h]
    \centering
    \begin{subfigure}[t]{0.48\textwidth}
         \centering
         \includegraphics[page=1]{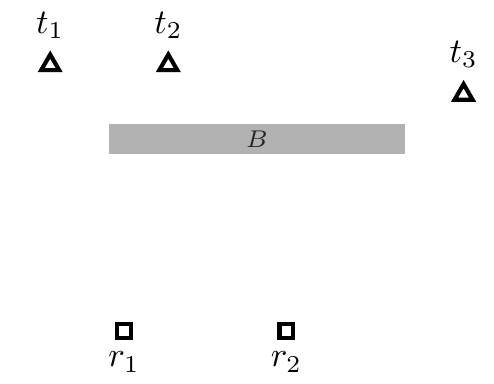}
        \caption{MRTA instance $\mathcal{M} = (R,T,L,p,c)$ with obstacle $B$.}
        \label{subfig:Example_InitialProb}
     \end{subfigure}
     \hfill
     \begin{subfigure}[t]{0.48\textwidth}
         \centering
         \includegraphics[page=2]{Figures/Eg2_H.pdf}
        \caption{Plan $\mathcal{P}$ output by $\textsc{Assign}(\mathcal{M})$, with $\edit{C}(\mathcal{P}) = 23.34$.}
        \label{subfig:Example_InitialSoln}
     \end{subfigure}
\caption{Initial MRTA problem and its solution. The positions of the tasks and robots are $p(r_1)=(2.5,0)$, $p(r_2) = (8,0)$, $p(t_1) = (0,9)$, $p(t_2) = (4,9)$ and $p(t_3) = (14,8)$.
}
\end{figure}

\begin{table*}[t] 
    \centering
    \caption{The interval family $\{(-\underline{\Delta}(e), \overline{\Delta}(e)]: e\in E(G(\mathcal{M}))\}$ output by $\textsc{AuctionSensitivity}(\mathcal{M}, \textsc{Auction}(\mathcal{M}))$.}
    \begin{tabularx}{\textwidth}{|l|| *{9}{>{\centering\arraybackslash} X|} }
    \hline
        Edge $e$ &
        $r_1t_1$ & $r_1t_2$ & $r_1t_3$ &
        $r_2t_1$ & $r_2t_2$ & $r_2t_3$ &
        $t_1t_2$ & $t_1t_3$ & $t_2t_3$
    \\
    \hline
        Cost $c(e)$ &
        9.34 & 9.85 & 14.06 &
        12.09 & 12.31 & 10 &
        4 & 14.04 & 10.05
    \\
        Maximum decrease $\underline{\Delta}(e)$ &
        9.34 & 0.25  & 4.04 &
        2.50 & 2.72  & 0.41 &
        4.00 & 4.01  & 0.02
    \\
        Maximum increase $\overline{\Delta}(e)$ &
          0.25   & $\infty$ & $\infty$ &
        $\infty$ & $\infty$ &   0.02   &
          5.60   & $\infty$ & $\infty$
    \\
    \hline
    \end{tabularx}
    \label{table:example_sensitivity}
\end{table*}

Suppose that during the execution of $\mathcal{P}$ we discover that our original data was inaccurate, and the obstacle is in a different position. Figure \ref{fig:example_updated} shows three different positions for the obstacle and the subsequent plans output. In Figure \ref{subfig:Example_reassign}, the obstacle is further east, in Figure \ref{subfig:Example_reorder}, it is further west; Figure \ref{subfig:Example_repeat} again shows the initial data for comparison.

\begin{figure}[hb]
    \centering
    \begin{subfigure}[t]{0.48\textwidth}
        \centering
        \includegraphics[page=3]{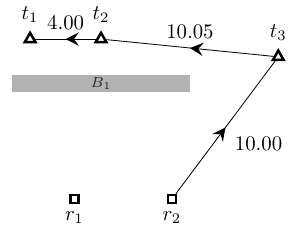}
        \caption{Instance $\mathcal{M}_1$, with obstacle $B_1$. Plan $\mathcal{P}_1 = \textsc{Assign}(\mathcal{M}_1)$ reassigns all tasks to $r_2$.}
        \label{subfig:Example_reassign}
     \end{subfigure}
     \hfill
         \begin{subfigure}[t]{0.48\textwidth}
        \centering
        \includegraphics[page=2]{Figures/Eg2_H.pdf}
        \caption{Instance $\mathcal{M}$ with obstacle $B$ and plan $\mathcal{P}$.}
        \label{subfig:Example_repeat}
     \end{subfigure}
     \hfill
     \begin{subfigure}[t]{0.48\textwidth}
        \centering
        \includegraphics[page=4]{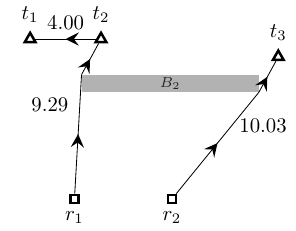}
        \caption{MRTA instance $\mathcal{M}_2$, with obstacle $B_2$. Plan $\mathcal{P}_2 = \textsc{Assign}(\mathcal{M}_2)$ reorders the tasks assigned to $r_1$.}
        \label{subfig:Example_reorder}
     \end{subfigure}
\caption{MRTA instances with different obstacle locations.}
\label{fig:example_updated}
\end{figure}

In both instances $\mathcal{M}_1$ and $\mathcal{M}_2$, our sensitivity analysis from Table \ref{table:example_sensitivity} tells us that this updated landscape data should trigger a replan:
for $\mathcal{M}_1$, we have 
$c_1(r_1t_1) = 10.18 > 9.34 + 0.25 = c(r_1t_1) + \overline{\Delta}(r_1t_1)$,
whereas for $\mathcal{M}_2$ we have 
$c_2(r_1t_2) = 9.29 < c(r_1t_2) - \underline{\Delta}(r_1t_2)$.\hfill{$\triangle$}

\end{exmp}

\section{Further work} \label{sec:FurtherWork}

In this paper, we provided a method to analyse the sensitivity of \textsc{Assign} to 
changes in the cost function on multiple edges. Specifically, we described how much the cost function can change before the list of winning edges $W$ output by the subroutine \textsc{Auction} (which forms the first step of \textsc{Assign}) also changes. See Theorem \ref{thm:LexBestFamily}. Notably, our sensitivity analysis does not consider the subroutine \textsc{DFShortcut}, which forms the second step of \textsc{Assign}.

A strength of this approach is that \textsc{DFShortcut} can be replaced by any alternative double-tree shortcutting algorithm on the same inputs, and our sensitivity analysis would still apply.
However, a weakness of this approach is that the intervals calculated in Theorem \ref{thm:LexBestFamily} are overly conservative for \textsc{Assign}, despite being tight for  \textsc{Auction}.

Thus there are two main ways in which our sensitivity analysis for \textsc{Assign} could be improved. 
First, observe that 
if two 
lists $W$ and $W'$ contain the \emph{same set} of winning edges in different orders, then our sensitivity analysis considers the output of \textsc{Auction} (i.e. $W$ or $W'$ respectively) to be different; whereas the two corresponding plans output by \textsc{DFShortcut} (and thus also by \textsc{Assign}) might be the same. 
Second, it's possible that two \emph{different sets} of winning edges 
could lead to $\textsc{DFShortcut}$ (and thus also \textsc{Assign}) outputting the same plan. Again, our current sensitivity analysis would flag these solutions as different.
%

Addressing the first of these extensions requires analysing the sensitivity of the \emph{set} of winning edges $\{w_k: w_k\in W\}$, rather than the \emph{list} $W=(w_k)$ to changes in the cost function. In the proof of Lemma \ref{lem:DoubleTreeCost}, we saw that the set 
of winning edges for a robot-task graph $G$ corresponds to a minimum spanning tree (MST) in an auxiliary graph $G'$. Thus this first proposed extension is equivalent to 
analysing the sensitivity of MST to multiple changes in the cost function. 
In 1982, Tarjan \cite{T_SensitivityAnalysisOfMST} showed how to analyse the sensitivity of MST to the cost of a single edge changing. 
Although this single edge variant has continued to receive attention
\cite{T_SAOfMST_2_worse, 
P_ImprovedSensitivityAnalysisMST}, 
no-one has yet succeeded in adapting Tarjan's result to multiple edges. 
%
%

\edit{An alternative research direction is to develop a similar sensitivity analysis for other algorithms. Possible candidates include: (1) algorithms for \emph{different MRTA problems} (such as \textsc{MinMax} or $\textsc{MinAvg}$, see \cite{L+_AuctionBasedMRR}), (2) different algorithms for the \emph{same MRTA problem} (\textsc{MinSum}), or (3) algorithms for subproblems to \textsc{MinSum}, such as Christofides's  $\frac{3}{2}$-approximation algorithm for TSP \cite{Christofides}.
}


\bibliographystyle{plain}
\bibliography{MRTA_bibliography.bib}

\end{document}